\documentclass{article}
\title{On the Notion of a Ribbon Quasi-Hopf Algebra}
\author{Yorck Sommerh\"auser}
\date{\it{To Susan Montgomery on the occasion of her 65th birthday}}

\usepackage{amsmath}
\usepackage{theorem}
\usepackage{amssymb}
\usepackage{tabularx}

\makeatletter
\renewcommand{\subsection}{\@startsection{subsection}{2}{0em}%
{\baselineskip}{-0em}{\bfseries\normalsize}}
\makeatother

\theoremstyle{plain}
\theorembodyfont{\rmfamily}

\newtheorem{thm}{Theorem}

\newtheorem{prop}[thm]{Proposition}
\newtheorem{lemma}[thm]{Lemma}

\newtheorem{pf}{Proof.}

\newtheorem{defn}[thm]{Definition}

\theoremstyle{break}
\theorembodyfont{\rmfamily}

\newcommand{\qed}{$\Box$}

\newcommand{\id}{\operatorname{id}}

\newcommand{\op}{\scriptstyle \operatorname{op}}
\newcommand{\cop}{\scriptstyle \operatorname{cop}}

\def\1{{(1)}}
\def\2{{(2)}}
\def\3{{(3)}}
\def\4{{(4)}}
\def\5{{(5)}}
\def\6{{(6)}}
\def\7{{(7)}}
\def\8{{(8)}}
\def\9{{(9)}}

\def\o{\otimes}

\def\da{\Delta}
\def\ea{\varepsilon}
\def\sa{S}
\def\A{1}

\begin{document}

\maketitle

\begin{abstract}
\hspace{-5mm}We show that two competing definitions of a ribbon quasi-Hopf algebra are actually equivalent. Along the way, we look at the Drinfel'd element from a new perspective and use this viewpoint to derive its fundamental properties.
\end{abstract}
\thispagestyle{empty}


\section*{Introduction} \label{Sec:Introd}
\addcontentsline{toc}{section}{Introduction}
While quasi-Hopf algebras were introduced by V.~G.~Drinfel'd (cf.~\cite{DrinfQuasiHopf}), the first authors to contemplate the notion of a ribbon quasi-Hopf algebra were D.~Altsch\"uler and A.~Coste (cf.~\cite{AltCoste}, Par.~4.1, p.~89). They define them as quasitriangular quasi-Hopf algebras with an additional central element, the ribbon element, that is subject to four axioms. However, as the authors point out themselves, these axioms are not completely satisfactory, as they neither reduce directly to the axioms of a ribbon Hopf algebra, in the case where the quasi-Hopf algebra happens to be an ordinary Hopf algebra, nor are in complete analogy to the axioms for a ribbon category. They therefore
analyzed their notion further and explained that, in the case
where the evaluation element~$\alpha$ is invertible, their axioms
are equivalent to a set of four different axioms which are considerably
closer to the notion of a ribbon Hopf algebra and the notion of 
a ribbon category. 

However, in the case of ribbon Hopf algebras, one of the four axioms is actually a consequence of the remaining axioms. Therefore D.~Bulacu, F.~Panaite, and F.~van Oystaeyen proposed a different definition of a ribbon quasi-Hopf algebra, leaving out this supposedly superfluous axiom
(cf.~\cite{BulPanOystTrace}, Def.~2.3, p.~6106). Again in the case where the evaluation element is invertible, they showed that this axiom really was superfluous, so that their definition was equivalent to the revised version of Altsch\"uler and Coste (cf.~\cite{BulPanOystTrace}, Prop.~5.5, p.~6119).

Of course, this raised the question whether the assumption on the invertibility of the evaluation element is really necessary to establish
these two equivalences, or whether this assumption was only made to simplify the argument. In the case of the first equivalence, between the two versions of the definition already proposed by Altsch\"uler and Coste, 
this question was addressed by D.~Bulacu and E.~Nauwelaerts, who showed that the assumption is not necessary (cf.~\cite{BulNauRib}, Thm.~3.1, p.~667). In a recent article, when using ribbon quasi-Hopf algebras to exemplify certain properties of modular data, the authors have claimed that this assumption is also not necessary for the second equivalence between the definition of Altsch\"uler and Coste and the definition of Bulacu, Panaite, and van Oystaeyen (cf.~\cite{SoZhuCent}, Cor.~5.1, p.~50). The purpose of the present article is to prove this claim.

To do this, we take a certain viewpoint, which is suitable not only for this proof, but also for similar questions: The R-matrix can be viewed
as a twist that takes the coproduct into the coopposite coproduct. However, while twisting leaves the antipode unchanged, the
coopposite coproduct naturally comes endowed with the inverse antipode.
The so-called Drinfel'd element now appears as the element that connects these two choices for the antipode of the coopposite quasi-Hopf algebra.
Viewing the Drinfel'd element in this way enables us not only to give
a relatively easy proof of our claim, but also allows us to give a new derivation of the fundamental properties of the Drinfel'd element in a comparatively short and conceptual way.

The article consists of two sections. The first, preliminary section contains a brief summary of the basic facts about quasi-bialgebras,
quasi-Hopf algebras, quasitriangularity, and twisting. However, we trace
more precisely than the available references how some elements already introduced in Drinfel'd's original article transform under twisting and 
other modifications, as this turns out to be crucial for our treatment.

The second section contains our main result, Theorem~\ref{RibQuasiHopf}. As explained above, we prove it by viewing the R-matrix as a twist, a viewpoint developed in Paragraph~\ref{RTwist}. The new proof of the fundamental properties of the Drinfel'd element also mentioned above is given in Paragraph~\ref{PropDrinf}. The article concludes with Proposition~\ref{OpCoop}, a formula for the image of the Drinfel'd element under the antipode. Although this formula was needed in our earlier proofs of Theorem~\ref{RibQuasiHopf}, it is not needed
in the proof presented here. We include it nonetheless, because it is
of independent interest and its proof nicely illustrates the ideas that we have developed.

In the following, we work over a base field that is denoted by~$K$.
All vector spaces that we will consider will be defined over this base field~$K$, and all tensor products will be taken over~$K$. With respect to enumeration, we use the convention that propositions, definitions, and similar items are referenced by the paragraph in which they occur; an 
additional third digit indicates a part of the corresponding item.
For example, a reference to Proposition~\ref{PropDrinf}.3 refers to the third assertion of the unique proposition in Paragraph~\ref{PropDrinf}.

\section{Preliminaries} \label{Sec:Prelim}
\subsection[Quasi-bialgebras]{} \label{QuasiBi}
Recall that a quasi-bialgebra is a quadruple
$(A,\da,\ea,\Phi)$, where~$A$ is an associative algebra
over our base field~$K$, whose multiplication and unit element we
have not explicitly listed as part of the structure elements. 
Out of the structure elements that we have listed explicitly,
two are algebra homomorphisms, namely~$\da: A \rightarrow A \o A$,
which we call the coproduct, and $\ea: A \rightarrow K$,  which we call the counit. The remaining structure element is the associator 
$\Phi \in  A \o A \o A$. These structure elements are required
to satisfy several axioms: Besides that $\Phi$ is required to be invertible, four equations have to be satisfied, which we now list. We give each equation a name that we will use in later references:
\begin{enumerate}
\item 
Quasi-coassociativity:
$(\id \o \da)\da(a) \Phi = \Phi (\da \o \id)\da(a)$

\item
Pentagon axiom:
\begin{align*}
(\id \o \id \o \da)(\Phi) (\da \o \id \o \id)(\Phi) 
= (\A \o \Phi) (\id \o \da \o \id)(\Phi) (\Phi \o \A)
\end{align*}

\item
Counitality:
$(\ea \o \id)\da(a) = a = (\id \o \; \ea)\da(a)$

\item 
Counit-associator axiom: $(\id \o \; \ea \o \id)(\Phi) = \A \o \A$
\end{enumerate}

Here, the first and the third equation are required for all $a \in A$. These axioms imply another property, which we call the 
counit-associator property:
\begin{prop}
$(\ea \o \id \o \id)(\Phi) = (\id \o \id \o \; \ea)(\Phi) = \A \o \A$
\end{prop}
\begin{pf}
This is proved in \cite{DrinfQuasiHopf}, Remark on p.~1422.
\qed
\end{pf}

We will use the version $\da(a) = a_\1 \o a_\2$ of the Heyneman-Sweedler sigma notation for the coproduct, and the notation
$\da^{\cop}(a) = a_\2 \o a_\1$ for the coopposite coproduct. Also, it will frequently be necessary to write~$\Phi$ and its inverse as a sum of decomposable tensors, which we do in the form
$$\Phi = \sum_{i=1}^n X_i \o Y_i \o Z_i \qquad  \qquad  
\Phi^{-1} = \sum_{j=1}^m \bar{X}_j \o \bar{Y}_j \o \bar{Z}_j$$
Because the number of decomposable tensors in these sums is never important in the sequel, we will also write such equations in slightly abbreviated forms, like $\Phi = \sum_i X_i \o Y_i \o Z_i$.

\subsection[Quasi-Hopf algebras]{} \label{QuasiHopf}
A quasi-bialgebra is a quasi-Hopf algebra if it is endowed with
three additional structure elements: An algebra anti-automorphism
$\sa: A \rightarrow A$, called the antipode, an element $\alpha \in A$,
called the evaluation element, and an element~$\beta \in A$, called the coevaluation element. The axioms that these structure elements have to satisfy are the following:
\begin{enumerate}
\item
Left antipode equation:
$\sa(a_\1) \alpha a_\2 = \ea(a) \alpha$
\item
Right antipode equation:
$a_\1 \beta \sa(a_\2) = \ea(a) \beta$
\item
Duality axiom: 
$\sum_i X_i \beta \sa(Y_i) \alpha Z_i = 1 = 
\sum_j \sa(\bar{X}_j) \alpha \bar{Y}_j \beta \sa(\bar{Z}_j)$
\end{enumerate}
These structure elements are compatible with the counit as follows:
\begin{lemma}
We have $\varepsilon(S(a)) = \varepsilon(a)$ and
$\varepsilon(\alpha) \varepsilon(\beta) = 1$.
\end{lemma}
\begin{pf}
The first assertion is proved in \cite{DrinfQuasiHopf}, Rem.~7, p.~1425.
The second follows by applying the counit to the
duality axiom.
\qed
\end{pf}
The antipode is also compatible with the coproduct and the associator. To formulate these compatibilities, we need to define two elements~$\gamma$ and~$\delta$ in the second tensor power of~$A$, which are in a sense analogues of the evaluation element~$\alpha$ and the coevaluation element~$\beta$:
\begin{align*}
\gamma &:= 
\sum_{i,j} \sa(\bar{X}_i Y_j) \alpha \bar{Y}_i Z_{j\1}\o 
\sa(X_j) \alpha \bar{Z}_i Z_{j\2} \\
\delta &:= 
\sum_{i,j} X_{i\1} \bar{X}_j \beta \sa(Z_i) \o X_{i\2} \bar{Y}_j \beta \sa(Y_i \bar{Z}_j) 
\end{align*}
From these elements, we derive the element
$$F := \sum_i (\sa(\bar{X}_{i\2}) \o \sa(\bar{X}_{i\1})) \gamma
\Delta(\bar{Y}_i \beta \sa(\bar{Z}_i))$$
which appears in the compatibility conditions in the following way:
\begin{prop}
$F$ is invertible with inverse
$$F^{-1} = \sum_i \Delta(\sa(\bar{X}_i) \alpha \bar{Y}_i) \delta (\sa(\bar{Z}_{i\2}) \o \sa(\bar{Z}_{i\1})) $$
and we have $\gamma = F \da(\alpha)$ and $\delta = \Delta(\beta) F^{-1}$.
The antipode is compatible with the coproduct via
$$\da(\sa(a)) = F^{-1} (\sa(a_\2) \o \sa(a_\1)) F$$
and with the associator via
$$\sum_i \sa(Z_i) \o \sa(Y_i) \o \sa(X_i) = 
(\A \o F) (\id \o \da)(F) \Phi (\da \o \id)(F^{-1}) (F^{-1} \o \A)$$
\end{prop}
\begin{pf}
This is proved in \cite{DrinfQuasiHopf}, Prop.~1.2, p.~1426. We note that
it is also shown there that the three properties
$\da(\sa(a)) = F^{-1} (\sa(a_\2) \o \sa(a_\1)) F$,
$\gamma = F \da(\alpha)$, and $\delta = \Delta(\beta) F^{-1}$
characterize~$F$ uniquely; even stronger, it suffices to check one
of the two conditions $\gamma = F \da(\alpha)$ and 
$\delta = \Delta(\beta) F^{-1}$.
\qed
\end{pf}

\subsection[Modifying the antipode]{} \label{ModAnti}
The antipode of a quasi-Hopf algebra is in general not unique;
it can be modified with the help of an invertible element~$x \in A$
by defining
$$\sa_x(a) := x S(a) x^{-1} \qquad \alpha_x := x \alpha \qquad 
\beta_x := \beta x^{-1}$$
It is easy to check that~$\sa_x$ is again an antipode for~$A$ with evaluation element~$\alpha_x$ and coevaluation element~$\beta_x$. However, this is
the only possible modification: If~$\sa'$ is an arbitrary new antipode for the quasi-Hopf algebra~$A$, with evaluation element~$\alpha'$ and coevaluation element~$\beta'$, then the element
$$x:= \sum_i \sa'(\bar{X}_i) \alpha' \bar{Y}_i \beta \sa(\bar{Z}_i)$$
is invertible with inverse 
$x^{-1} = \sum_i \sa(\bar{X}_i) \alpha \bar{Y}_i \beta' \sa'(\bar{Z}_i)$,
and we have $\sa' =\sa_x$, $\alpha' = \alpha_x$, and $\beta' = \beta_x$.
This fact, which will be important in the sequel, is proved in 
\cite{DrinfQuasiHopf}, Prop.~1.1, p.~1425.

By modifying the antipode as indicated by an invertible element~$x$, we of course indirectly modify all other elements derived from it; in particular
the elements~$\gamma$, $\delta$, and~$F$ introduced in Paragraph~\ref{QuasiHopf}. The modified elements, which we denote by~$\gamma_x$, $\delta_x$, and~$F_x$, can be expressed in terms of the
unmodified elements as follows:
\begin{prop}
$$\gamma_x  = (x \o x) \gamma \qquad \qquad
\delta_x  = \delta (x^{-1} \o x^{-1}) \qquad \qquad
F_x  = (x \o x) F \Delta(x^{-1})$$
\end{prop}
\begin{pf}
The form of~$\gamma_x$ follows directly from the definition:
\begin{align*}
\gamma_x = \sum_{i,j} 
x S(\bar{X}_i Y_j) x^{-1} (x \alpha) \bar{Y}_i Z_{j\1} \o x S(X_j) x^{-1} (x \alpha) \bar{Z}_i Z_{j\2}
= (x \o x) \gamma
\end{align*}
Similarly, the definition of~$\delta_x$ is
\begin{align*}
\delta_x = \sum_{i,j} 
X_{i\1} \bar{X}_j (\beta x^{-1}) x S(Z_{i}) x^{-1} \o 
X_{i\2} \bar{Y}_j (\beta x^{-1}) x S(Y_i \bar{Z}_j) x^{-1}
\end{align*}
which immediately yields the second assertion. Finally, since
\begin{align*}
F_x &= 
\sum_i (\sa_x(\bar{X}_{i\2}) \o \sa_x(\bar{X}_{i\1})) \gamma_x
\Delta(\bar{Y}_i \beta_x \sa_x(\bar{Z}_i)) \\
&= \sum_i (x \sa(\bar{X}_{i\2}) x^{-1} \o x \sa(\bar{X}_{i\1}) x^{-1}) 
(x \o x) \gamma
\Delta(\bar{Y}_i (\beta x^{-1}) x \sa(\bar{Z}_i) x^{-1}) \\
&= (x \o x) \sum_i (\sa(\bar{X}_{i\2}) \o \sa(\bar{X}_{i\1})) 
\gamma \Delta(\bar{Y}_i \beta \sa(\bar{Z}_i) x^{-1})
= (x \o x) F \Delta(x^{-1})
\end{align*}
the third assertion also holds.
\qed
\end{pf}

\subsection[The coopposite quasi-Hopf algebra]{} \label{Coop}
With every quasi-Hopf algebra~$A$, one can associate another 
quasi-Hopf algebra~$A^{\cop}$, which has the same product as~$A$,
but the coopposite coproduct. For this quasi-Hopf algebra, the counit is unchanged, the associator is changed to
$\sum_i \bar{Z}_i \o \bar{Y}_i \o \bar{X}_i$, the antipode is changed
to its inverse~$\sa^{-1}$, the evaluation element is changed 
to~$\sa^{-1}(\alpha)$, and the coevaluation element is changed to~$\sa^{-1}(\beta)$ (cf.~\cite{DrinfQuasiHopf}, Rem.~4, p.~1424; \cite{Kas}, Exerc.~XV.6.2, p.~381).

As in Paragraph~\ref{ModAnti}, this modification of the defining structure elements also leads to a modification of the elements~$\gamma$, $\delta$, and~$F$. In this case, however, we do not introduce a special notation for the new elements formed in~$A^{\cop}$, because their relation
to the original elements is so simple: The new elements are 
$(S^{-1} \o S^{-1})(\gamma)$, $(S^{-1} \o S^{-1})(\delta)$, and 
$(S^{-1} \o S^{-1})(F)$. To see this in the case of~$\gamma$, we use an alternative description of~$\gamma$ given in \cite{DrinfQuasiHopf}, Lem.~1, p.~1427, which yields
\begin{align*}
(S^{-1} \o S^{-1})(\gamma) &= 
(S^{-1} \o S^{-1})(\sum_{i,j} 
S(Y_i \bar{X}_{j\2}) \alpha Z_i \bar{Y}_j \o 
S(X_i \bar{X}_{j\1}) \alpha \bar{Z}_j) \\
&= \sum_{i,j} 
S^{-1}(Z_i \bar{Y}_j) S^{-1}(\alpha) Y_i \bar{X}_{j\2} \o S^{-1}(\bar{Z}_j) S^{-1}(\alpha) X_i \bar{X}_{j\1}
\end{align*}
But this last term is just what we get if we form~$\gamma$ in~$A^{\cop}$
according to the original definition in Paragraph~\ref{QuasiHopf}. 

In the case of~$\delta$, we argue similarly: An alternative formula
given in \cite{DrinfQuasiHopf}, loc.~cit.\ implies that 
\begin{align*}
(S^{-1} \o S^{-1})(\delta) &= 
(S^{-1} \o S^{-1})(\sum_{i,j} \bar{X}_{i} \beta S(\bar{Z}_{i\2} Z_j) \o 
\bar{Y}_i X_j\beta S(\bar{Z}_{i\1} Y_j)) \\
&= \sum_{i,j} 
\bar{Z}_{i\2} Z_j S^{-1}(\beta) S^{-1}(\bar{X}_{i}) \o 
\bar{Z}_{i\1} Y_j S^{-1}(\beta) S^{-1}(\bar{Y}_i X_j)
\end{align*}
which is again what we get if we form~$\delta$ in~$A^{\cop}$
according to the original definition in Paragraph~\ref{QuasiHopf}. 

In the case of~$F$, we argue differently: If we apply $S^{-1} \o S^{-1}$ to the equation $(S \o S)(\Delta^{\cop}(a)) = F \Delta(S(a)) F^{-1}$ in Proposition~\ref{QuasiHopf} and replace~$a$ by~$S^{-1}(a)$, we get
$$\Delta^{\cop}(S^{-1}(a)) = 
(S^{-1} \o S^{-1})(F^{-1}) (S^{-1} \o S^{-1}) (\Delta(a))
(S^{-1} \o S^{-1})(F) $$
Similarly, if we apply $S^{-1} \o S^{-1}$ to the equation 
$\gamma = F \Delta(\alpha)$ in the same proposition and use what we have just established, we get 
\begin{align*}
(S^{-1} \o S^{-1})(\gamma) &= 
(S^{-1} \o S^{-1})(\Delta(\alpha)) (S^{-1} \o S^{-1})(F) \\
&= (S^{-1} \o S^{-1})(F) \Delta^{\cop}(S^{-1}(\alpha))
\end{align*}
Finally, if we treat the equation $\delta = \Delta(\beta) F^{-1}$
in the same way, we get
$$(S^{-1} \o S^{-1})(\delta) = 
\Delta^{\cop}(S^{-1}(\beta)) (S^{-1} \o S^{-1})(F^{-1})$$
But this establishes our assertion, since it shows that 
$(S^{-1} \o S^{-1})(F)$ has the characteristic properties of the element~$F$ in~$A^{\cop}$, as described in Paragraph~\ref{QuasiHopf}.

\subsection[Quasitriangularity]{} \label{QuasiTri}
A quasi-Hopf algebra is called quasitriangular if it is endowed with a
so-called R-matrix, which is an invertible element
$R = \sum_l s_l \o t_l \in A \o A$ that satisfies the following three
conditions:
\begin{enumerate}
\item 
Quasi-cocommutativity: $\da^{\cop}(a) R = R \da(a)$

\item
Left hexagon axiom:
\begin{align*}
(\da \o \id)(R) = \sum_{i,j,k,l,q}  
Y_i s_l \bar{X}_j X_k \o Z_i \bar{Z}_j s_q Y_k \o 
X_i t_l \bar{Y}_j t_q Z_k 
\end{align*}

\item
Right hexagon axiom:
\begin{align*}
&(\id \o \da)(R) 
= \sum_{i,j,k,l,q} \bar{Z}_i s_l Y_j s_q \bar{X}_k 
\o \bar{X}_i X_j t_q \bar{Y}_k \o \bar{Y}_i t_l Z_j \bar{Z}_k 
\end{align*}
\end{enumerate}

Note that the right-hand side in the hexagon axioms factors
completely; for example, the right-hand side in the left hexagon
axiom is the product of the tensors
$\sum_i Y_i \o Z_i \o X_i$, $\sum_l s_l \o \A \o t_l$,
$\sum_j \bar{X}_j \o \bar{Z}_j \o \bar{Y}_j$, 
$\sum_q \A \o s_q \o t_q$, and $\sum_k X_k \o Y_k \o Z_k$.

The hexagon axioms obviously constitute a compatibility condition
between the R-matrix and the coproduct. But the R-matrix is also compatible
with the counit and the antipode: Denoting by~$F'$ the image of~$F$ under the interchange of the two tensor factors, we have
\begin{lemma}
$$(\ea \o \id)(R) = \A \qquad \qquad (\id \o \, \ea)(R) = \A
\qquad \qquad (\sa \o \sa)(R) = F'RF^{-1}$$
\end{lemma}
\begin{pf}
The equations involving the counit are proved in \cite{DrinfQuasiHopf}, Rem.~2, p.~1440; they are also stated in \cite{AltCoste}, Eq.~(2.23), p.~87. The equation involving the antipode was stated in \cite{AltCoste}, Eq.~(4.22), p.~96 and proved in \cite{HN}, Cor.~2.2, p.~559. A proof
without the graphical calculus was given in \cite{BulNauRib}, Lem.~2.3, p.~663. These references also list additional compatibility conditions between the R-matrix and the antipode.
\qed
\end{pf}

From the R-matrix, we derive a special element~$u$, called the 
Drinfel'd element. It is defined as 
$$u := \sum_{i,l} 
\sa(\bar{Y}_i \beta \sa(\bar{Z}_i)) \sa(t_l) \alpha s_l \bar{X}_i$$
(cf.~\cite{AltCoste}, Eq.~(3.2), p.~87; \cite{Kas}, Exerc.~XV.6.5, p.~381). 
This ad hoc definition may appear unmotivated at this point; we will
put it in its context in Paragraph~\ref{RTwist}. Although we could set
down the fundamental properties of the Drinfel'd element here, as they appear in literature, we defer this to Paragraph~\ref{PropDrinf}, where we will actually reconfirm them from the viewpoint developed in Paragraph~\ref{RTwist}, as this viewpoint allows for a proof that is in our opinion shorter and more conceptual. Here we only record how the
Drinfel'd element changes if the antipode is modified by an invertible
element~$x$ as explained in Paragraph~\ref{ModAnti}. The new Drinfel'd element~$u_x$ relates to the old Drinfel'd element~$u$ as follows:
\begin{prop}
$u_x = x \sa(x^{-1}) u$
\end{prop}
\begin{pf}
As we have
\begin{align*}
u_x &= \sum_{i,l} 
\sa_x(\bar{Y}_i \beta_x \sa_x(\bar{Z}_i)) \sa_x(t_l) \alpha_x s_l \bar{X}_i \\
&= \sum_{i,l} 
x \sa(\bar{Y}_i (\beta x^{-1}) (x \sa(\bar{Z}_i) x^{-1})) x^{-1} 
(x \sa(t_l) x^{-1}) (x \alpha) s_l \bar{X}_i \\
&= \sum_{i,l} 
x \sa(\bar{Y}_i \beta \sa(\bar{Z}_i) x^{-1}) \sa(t_l) \alpha s_l \bar{X}_i
= x \sa(x^{-1}) u
\end{align*}
we see that this follows directly from the definition.
\qed
\end{pf}

\subsection[Twisting]{} \label{Twist}
Quasi-Hopf algebras can be twisted to generate new quasi-Hopf algebras.
The ingredient that we need for this is a twisting element; i.e., an invertible element $T \in A \o A$ in the second tensor power of our quasi-Hopf algebra~$A$ that satisfies the condition 
$(\ea \o \id)(T) = (\id \o \; \ea)(T) = \A$.
If we then introduce the new coproduct
$$\Delta_T(a) := T \da(a) T^{-1}$$
and the new associator
$$\Phi_T := 
(\A \o T)(\id \o \da)(T) \Phi (\da \o \id)(T^{-1})(T^{-1} \o \A)$$
but leave the counit and the antipode unchanged, we get again a quasi-Hopf algebra, at least if we introduce a new evaluation
element~$\alpha_T$ and a new coevaluation element~$\beta_T$ via
$$\alpha_T := \sum_i \sa(\bar{f}_i) \alpha \bar{g}_i \qquad \qquad
\beta_T := \sum_i f_i \beta \sa(g_i)$$
where we have used the notation $T = \sum_i f_i \o g_i$
and $T^{-1} = \sum_i \bar{f}_i \o \bar{g}_i$ (cf.~\cite{DrinfQuasiHopf}, Rem.~5, p.~1425; \cite{Kas}, Exerc.~XV.6.4, p.~381).

As a consequence of these modifications, we also get, according to our definitions in Paragraph~\ref{QuasiHopf}, new elements~$\gamma_T$, $\delta_T$, and~$F_T$. As we will show now, these new elements can be expressed in terms of the original elements~$\gamma$, $\delta$, and~$F$. If we denote, as for~$F$, by~$T'$ the image of~$T$ under the interchange of the two tensor factors, the corresponding expressions look, in a slightly implicit form, as follows:

\pagebreak

\begin{prop}
\enlargethispage{0.1cm}
\begin{enumerate}
\item[]
\item 
$\displaystyle
(\sa \o \sa)(T') \gamma_T T = 
\sum_{i} (\sa \o \sa)(\da^{\cop}(\bar{f}_{i}))
\gamma \da(\bar{g}_i)$

\item
$\displaystyle
T^{-1} \delta_T (\sa \o \sa)(T'^{-1})
= \sum_{i} \da(f_{i}) \delta (\sa \o \sa)(\da^{\cop}(g_i))$

\item
$\displaystyle
F_T = (\sa \o \sa)(T'^{-1}) F T^{-1}$
\end{enumerate}
\end{prop}
\begin{pf}
\newcounter{num}
\begin{list}{(\arabic{num})}{\usecounter{num} \leftmargin0cm \itemindent5pt}
\item
We use the Sweedler notation $\Delta_T(a) = a_{[1]} \o a_{[2]}$
for the twisted coproduct, and primes for the twisted associator; i.e., we write
$$\Phi_T = \sum_{i} X'_i \o Y'_i \o Z'_i \qquad \qquad
\Phi_T^{-1} = \sum_{j} \bar{X}'_j \o \bar{Y}'_j \o \bar{Z}'_j$$
With this notation, the definition of~$\gamma_T$ reads
\begin{align*}
\gamma_T &= \sum_{i,j} 
\sa(\bar{X}'_i Y'_j) \alpha_T \bar{Y}'_i Z'_{j[1]} \o \sa(X'_j) \alpha_T \bar{Z}'_i Z'_{j[2]} \\
&= \sum_{i,j,k,l} 
\sa(\bar{f}_k \bar{X}'_i Y'_j) \alpha \bar{g}_k \bar{Y}'_i Z'_{j[1]} \o \sa(\bar{f}_l X'_j) \alpha \bar{g}_l \bar{Z}'_i Z'_{j[2]}
\end{align*}
If we multiply this from the right by $T = \sum_q f_q \o g_q$ and use the fact that $\da_T(a) T = T \da(a)$, we get
\begin{align*}
\gamma_T T &= \sum_{i,j,k,l,q} 
\sa(\bar{f}_k \bar{X}'_i Y'_j) \alpha \bar{g}_k \bar{Y}'_i f_q Z'_{j\1} \o \sa(\bar{f}_l X'_j) \alpha \bar{g}_l \bar{Z}'_i g_q Z'_{j\2}
\end{align*}
But from the definition of the twisted associator, we have
\begin{align*}
\sum_{i,k,q} 
\bar{f}_k \bar{X}'_i \o \bar{g}_k \bar{Y}'_i f_q \o  \bar{Z}'_i g_q
&= (T^{-1} \o 1) \Phi_T^{-1}(1 \o T) \\
&= (\da \o \id)(T) \Phi^{-1} (\id \o \da)(T^{-1}) 
\end{align*}
If we insert this into our expression, the term
$(\da \o \id)(T)$ cancels, and we get
\begin{align*}
\gamma_T T &= \sum_{i,j,k,l} 
\sa(\bar{X}_i \bar{f}_k Y'_j) \alpha \bar{Y}_i \bar{g}_{k\1} Z'_{j\1} \o \sa(\bar{f}_l X'_j) \alpha \bar{g}_l \bar{Z}_i \bar{g}_{k\2}  Z'_{j\2} 
\end{align*}
Multiplying from the left by~$(\sa \o \sa)(T')$ yields
\begin{align*}
(\sa \o \sa)(&T') \gamma_T T = \\
&\sum_{i,j,k,l,q} 
\sa(\bar{X}_i \bar{f}_k Y'_j g_q) \alpha \bar{Y}_i \bar{g}_{k\1} Z'_{j\1} \o \sa(\bar{f}_l X'_j f_q) \alpha \bar{g}_l \bar{Z}_i \bar{g}_{k\2}  Z'_{j\2} 
\end{align*}
Now we have, again from the definition of the twisted associator, that
\begin{align*}
&\sum_{j,k,q} X'_j f_q \o \bar{f}_k Y'_j g_q \o \bar{g}_{k} Z'_{j}
= (1 \o T^{-1}) \Phi_T (T \o 1) \\
&= (\id \o \da)(T) \Phi (\da \o \id)(T^{-1}) 
= \sum_{j,k,q} f_{k} X_j \bar{f}_{q\1} 
\o g_{k\1} Y_j \bar{f}_{q\2} 
\o g_{k\2} Z_j \bar{g}_q 
\end{align*}
Inserting this, our expression becomes
\begin{align*}
(\sa \o \sa)(T') \gamma_T T = \sum_{i,j,k,l,q} 
\sa(\bar{X}_i g_{k\1} &Y_j \bar{f}_{q\2}) \alpha \bar{Y}_i g_{k\2\1} Z_{j\1} \bar{g}_{q\1} \\
& \o \sa(\bar{f}_l f_{k} X_j \bar{f}_{q\1}) \alpha \bar{g}_l \bar{Z}_i g_{k\2\2} Z_{j\2} \bar{g}_{q\2} 
\end{align*}
Using quasi-coassociativity, we can write this as
\begin{align*}
(\sa \o \sa)(T') \gamma_T T = \sum_{i,j,k,l,q} 
\sa(g_{k\1\1} \bar{X}_i & Y_j \bar{f}_{q\2}) \alpha  g_{k\1\2} \bar{Y}_i Z_{j\1} \bar{g}_{q\1} \\
&\o \sa(\bar{f}_l f_{k} X_j \bar{f}_{q\1}) \alpha \bar{g}_l g_{k\2} \bar{Z}_i Z_{j\2} \bar{g}_{q\2} 
\end{align*}
Here we can use the left antipode equation on the part
$\sa(g_{k\1\1}) \alpha  g_{k\1\2}$, and after that the summations
over~$k$ and~$l$ cancel, so that we are left with
\begin{align*}
(\sa \o \sa)(T') \gamma_T T &= \sum_{i,j,q} 
\sa(\bar{X}_i Y_j \bar{f}_{q\2}) \alpha \bar{Y}_i Z_{j\1} \bar{g}_{q\1}  \o \sa(X_j \bar{f}_{q\1}) \alpha \bar{Z}_i Z_{j\2} \bar{g}_{q\2} \\
&=  \sum_{q} (\sa(\bar{f}_{q\2}) \o \sa(\bar{f}_{q\1}))
\gamma (\bar{g}_{q\1}  \o \bar{g}_{q\2}) 
\end{align*}
which is the first assertion.

\item
The form of~$\delta_T$ can be established by a very similar computation.
However, this computation can be avoided by using the argument that we present now. In this approach, we redefine~$F_T$ to be what we claim it is
according to the third assertion, i.e., we redefine it as
$F_T := (\sa \o \sa)(T'^{-1}) F T^{-1}$. By Proposition~\ref{QuasiHopf}, the original element~$F$ satisfies 
$(\sa \o \sa)(\da^{\cop}(a)) = F \da(\sa(a)) F^{-1}$, so~$F_T$ satisfies
\begin{align*}
(\sa \o \sa)(&\Delta_T^{\cop}(a)) = (\sa \o \sa)(T' \da^{\cop}(a) T'^{-1}) \\ 
&= (\sa \o \sa)(T'^{-1}) F \da(\sa(a)) F^{-1} (\sa \o \sa)(T')= F_T \da_T(\sa(a)) F_T^{-1}
\end{align*}
Furthermore, from the first assertion and the properties of~$F$ we have that
\begin{align*}
(\sa \o \sa)(T') \gamma_T T &= 
\sum_{i} (\sa \o \sa)(\da^{\cop}(\bar{f}_{i}))
\gamma \da(\bar{g}_i) = 
\sum_{i} F \da(\sa(\bar{f}_{i})) F^{-1} \gamma \da(\bar{g}_i) \\
&= F \sum_{i} \da(\sa(\bar{f}_{i})) \da(\alpha) \da(\bar{g}_i)
= F \da(\alpha_T)
\end{align*}
so that 
$\gamma_T = (\sa \o \sa)(T'^{-1}) F \da(\alpha_T) T^{-1} 
= F_T \Delta_T(\alpha_T)$. 
But we know from Paragraph~\ref{QuasiHopf} that the two properties
that we have just established characterize~$F_T$, in other words,
the third assertion of our proposition holds. 

\item
But then the equation $\delta_T = \Delta_T(\beta_T) F_T^{-1}$ holds by Proposition~\ref{QuasiHopf}. Inserting the form of~$F_T$, this says that
$\delta_T = T \da(\beta_T) F^{-1} (\sa \o \sa)(T')$, so that the left-hand side of the second assertion of our proposition is
\begin{align*}
T^{-1} \delta_T (\sa \o \sa)(T'^{-1}) = \da(\beta_T) F^{-1}
= \sum_{i} \da(f_{i}) \da(\beta) \da(\sa(g_i)) F^{-1}
\end{align*}
Using the properties of~$F$ again, we can rewrite this in the 
form
\begin{align*}
T^{-1} \delta_T (\sa \o \sa)(T'^{-1}) 
&= \sum_{i} \da(f_{i}) \da(\beta) F^{-1} (\sa \o \sa)(\da^{\cop}(g_i)) \\
&= \sum_{i} \da(f_{i}) \delta (\sa \o \sa)(\da^{\cop}(g_i))
\end{align*}
where the last step uses the original equation $\delta = \da(\beta) F^{-1}$
from Proposition~\ref{QuasiHopf}. But this is exactly the second assertion
of our proposition.
\qed
\end{list}
\end{pf}

If~$A$ is quasitriangular, its twist is also quasitriangular, with respect to the new R-matrix $R_T := T' R T^{-1}$ (cf.~\cite{DrinfQuasiHopf}, Eq.~(3.11), p.~1439; \cite{Kas}, Prop.~XV.3.6, p.~376). This new R-matrix in principle also gives rise to a new Drinfel'd element~$u_T$. However, this new element coincides with the original one:
\begin{lemma}
$u_T = u$
\end{lemma}
\begin{pf}
This is proved in \cite{BulPanOystTrace}, Lem.~4.2, p.~6115. 
\qed
\end{pf}

\section{Ribbon quasi-Hopf algebras} \label{Sec:RibQuasiHopf}
\subsection[The R-matrix as a twist]{} \label{RTwist}
We can also relate quasitriangularity and twisting in another way:
As the twisting element~$T$, we can choose the R-matrix~$R$, because
Lemma~\ref{QuasiTri} asserts that the R-matrix satisfies
the conditions that a twist element should satisfy. By definition,
the twisted coproduct is just the coopposite coproduct, which we have
discussed in Paragraph~\ref{Coop}. However, not all of the other structure elements match: Although the Yang-Baxter equation (cf.~\cite{Kas}, Cor.~XV.2.3, p.~372) yields that the twisted associator is also
$\sum_i \bar{Z}_i \o \bar{Y}_i \o \bar{X}_i$, the antipode remains the same, and is not changed to its inverse, and for the evaluation 
element~$\hat{\alpha}$ and the coevaluation element~$\hat{\beta}$
we find the expressions 
$$\hat{\alpha} = \sum_l \sa(\bar{s}_l) \alpha \bar{t}_l \qquad  \qquad
\hat{\beta} = \sum_{l} s_l \beta \sa(t_l)$$
where we have, as before, used the notation 
$R = \sum_{l} s_l \o t_l$ and 
$R^{-1} = \sum_l \bar{s}_l \o \bar{t}_l$.

This is, however, not a contradiction; we have already discussed in
Paragraph~\ref{ModAnti} that the antipode of a quasi-Hopf algebra is not unique, and we have also explained there how the structures are related: The element
$$\hat{u} := 
\sum_i \sa(Z_i) \hat{\alpha} Y_i \sa^{-1}(\beta) \sa^{-1}(X_i)$$
is invertible with inverse 
$\hat{u}^{-1} = 
\sum_i \sa^{-1}(Z_i) \sa^{-1}(\alpha) Y_i \hat{\beta} \sa(X_i)$, and 
we have
$$\sa(a) = \hat{u} \sa^{-1}(a) \hat{u}^{-1} \qquad \qquad
\hat{\alpha} = \hat{u} \sa^{-1}(\alpha)\qquad \qquad 
\hat{\beta} = \sa^{-1}(\beta) \hat{u}^{-1}$$

Now we can associate with every R-matrix another one: It follows directly from the definition in Paragraph~\ref{QuasiTri} that~$R'^{-1}$ is also
an R-matrix for~$A$, where~$R'$ denotes, as for~$F$ and~$T$ before, the image of~$R$ under the interchange of the two tensor factors. We can therefore also use this R-matrix to twist the coproduct into the coopposite coproduct. In this case, the twisted associator is again $\sum_i \bar{Z}_i \o \bar{Y}_i \o \bar{X}_i$, the antipode remains unchanged, and for the evaluation 
element~$\check{\alpha}$ and the coevaluation element~$\check{\beta}$
we find the expressions 
$$\check{\alpha} = \sum_l \sa(t_l) \alpha s_l \qquad  \qquad
\check{\beta} = \sum_{l} \bar{t}_l \beta \sa(\bar{s}_l)$$
Also the discussion in Paragraph~\ref{ModAnti} applies again to tell us the relation of the structures: The element
$$\check{u} := 
\sum_i \sa(Z_i) \check{\alpha} Y_i \sa^{-1}(\beta) \sa^{-1}(X_i)$$
is invertible with inverse 
$\check{u}^{-1} = 
\sum_i \sa^{-1}(Z_i) \sa^{-1}(\alpha) Y_i \check{\beta} \sa(X_i)$, and 
we have
$$\sa(a) = \check{u} \sa^{-1}(a) \check{u}^{-1} \qquad \qquad
\check{\alpha} = \check{u} \sa^{-1}(\alpha) \qquad \qquad 
\check{\beta} = \sa^{-1}(\beta) \check{u}^{-1}$$

It is to be expected that there is a connection between these two ways of twisting. A first connection involves the evaluation and the coevaluation elements:
\begin{lemma}
For the evaluation elements, we have
$$\sa^{-1}(\check{\alpha}) = \hat{u}^{-1} \alpha \qquad \qquad \qquad
\sa^{-1}(\hat{\alpha}) = \check{u}^{-1} \alpha$$
For the coevaluation elements, we have
$$\sa^{-1}(\check{\beta}) = \beta \hat{u} \qquad \qquad \qquad
\sa^{-1}(\hat{\beta}) = \beta \check{u}$$
\end{lemma}
\begin{pf}
It is easy to solve the definitions of $\hat{\alpha}$ and $\check{\alpha}$
for~$\alpha$; we find
$$\alpha = \sum_l \sa(s_l) \hat{\alpha} t_l \qquad \qquad
\alpha = \sum_l \sa(\bar{t}_l) \check{\alpha} \bar{s}_l$$
If we apply the inverse antipode to the definition of~$\check{\alpha}$
and use the preceding formulas, we therefore get
\begin{align*}
\sa^{-1}(\check{\alpha}) = \sum_l \sa^{-1}(s_l) \sa^{-1}(\alpha) t_l = \sum_l \sa^{-1}(s_l) \hat{u}^{-1} \hat{\alpha} t_l
= \hat{u}^{-1} \sum_l \sa(s_l) \hat{\alpha} t_l= \hat{u}^{-1} \alpha
\end{align*}
The formula $\sa^{-1}(\hat{\alpha}) = \check{u}^{-1} \alpha$ can be
established by a similar computation, but on the other hand, it also
follows from the first equation by interchanging~$R$ and~$R'^{-1}$.

The coevaluation elements can be treated similarly: Solving their 
definitions for~$\beta$, we find
$$\beta = \sum_l \bar{s}_l \hat{\beta} \sa(\bar{t}_l) \qquad \qquad
\beta = \sum_l t_l \check{\beta} \sa(s_l)$$
If we apply the inverse antipode to the definition of~$\check{\beta}$
and use the preceding formulas, we therefore get
\begin{align*}
\sa^{-1}(\check{\beta}) = 
\sum_l \bar{s}_l \sa^{-1}(\beta) \sa^{-1}(\bar{t}_l) = 
\sum_l \bar{s}_l \hat{\beta} \hat{u} \sa^{-1}(\bar{t}_l) = 
\sum_l \bar{s}_l \hat{\beta} \sa(\bar{t}_l) \hat{u} = \beta \hat{u}
\end{align*}
Again, the formula $\sa^{-1}(\hat{\beta}) = \beta \check{u}$ can be
established by a similar computation, or viewed as a consequence by interchanging~$R$ and~$R'^{-1}$.
\qed
\end{pf}

There is also a direct connection between the elements~$\hat{u}$ and~$\check{u}$, and, what is important for us, there is a connection to the Drinfel'd element~$u$:
\begin{prop}
$u = \check{u} = \sa(\hat{u}^{-1})$
\end{prop}
\begin{pf}
From the preceding lemma, we get that
$$\hat{u}^{-1} \alpha = \sa^{-1}(\check{\alpha}) 
= \sa^{-2}(\alpha) \sa^{-1}(\check{u})
= \hat{u}^{-1} \alpha \hat{u} \sa^{-1}(\check{u})$$
so that $\alpha = \alpha \hat{u} \sa^{-1}(\check{u})$.
Now the square of the antipode is both conjugation with
$\hat{u}$ and conjugation with $\sa^{-1}(\check{u}^{-1})$, so
that $\hat{u} \sa^{-1}(\check{u})$ is a central element.
But then the duality axiom implies that
$$\hat{u} \sa^{-1}(\check{u}) = 
\sum_i X_i \beta \sa(Y_i) \alpha \hat{u} \sa^{-1}(\check{u}) Z_i =  
\sum_i X_i \beta \sa(Y_i) \alpha Z_i = 1$$
This shows that $\hat{u}^{-1} = \sa^{-1}(\check{u})$
and therefore $\check{u} = \sa(\hat{u}^{-1})$.

For the assertion about the Drinfel'd element, we first note that with our new terminology we can rewrite its definition, given in Paragraph~\ref{QuasiTri}, in the form
$$u = \sum_{i} 
\sa(\bar{Y}_i \beta \sa(\bar{Z}_i))  \check{\alpha} \bar{X}_i$$
Applying the inverse antipode and using that 
$\sa^{-1}(\check{\alpha}) = \hat{u}^{-1} \alpha$
by the preceding lemma, we get
\begin{align*}
\sa^{-1}(u) = 
\sum_i \sa^{-1}(\bar{X}_i) \hat{u}^{-1} \alpha 
\bar{Y}_i \beta \sa(\bar{Z}_i) = 
\hat{u}^{-1} \sum_i \sa(\bar{X}_i) \alpha 
\bar{Y}_i \beta \sa(\bar{Z}_i) = \hat{u}^{-1}
\end{align*}
where the last step follows from the duality axiom.
This shows that $u = \sa(\hat{u}^{-1})$, as asserted.
\qed
\end{pf}

We note that this proposition and the preceding lemma imply immediately that $\check{\alpha} = \sa(\alpha) u$, which is an identity that appears in
\cite{AltCoste}, Eq.~(3.9), p.~88.

\subsection[The properties of the Drinfel'd element]{} \label{PropDrinf}
The choice of the R-matrix~$R$ as the twisting element~$T$ does not only
lead to the elements~$\hat{\alpha}$, $\hat{\beta}$, and~$\hat{u}$, but also, as we saw in Paragraph~\ref{Twist}, to new versions of the
elements~$\gamma$, $\delta$, and~$F$, which we denote by
$\hat{\gamma}$, $\hat{\delta}$, and~$\hat{F}$. Similarly, the choice of~$R'^{-1}$ as the twisting element~$T$ leads to new versions of these elements that we denote by $\check{\gamma}$, $\check{\delta}$, and~$\check{F}$. We have seen in Proposition~\ref{Twist} how the new
elements can be expressed in terms of the old ones; we record here
only the form of~$\hat{F}$ and~$\check{F}$, where this proposition yields that
$$\hat{F} = (\sa \o \sa)(R'^{-1}) F R^{-1} \qquad \qquad 
\check{F} = (\sa \o \sa)(R) F R'$$
We now use all of this to derive the fundamental properties of the Drinfel'd element~$u$, as promised in the introduction and in Paragraph~\ref{QuasiTri}.
These fundamental properties are the following:
\begin{prop}
$u$ is invertible. Moreover, we have
\begin{enumerate}
\item
$\ea(u) = 1$

\item 
$\sa^2(a) = u a u^{-1}$

\item
$\da(u) = F^{-1} ((\sa \o \sa)(F')) (u \o u) (R'R)^{-1}$
\end{enumerate}
\end{prop}
\begin{pf}
By Proposition~\ref{RTwist}, we have $u=\check{u}$, and we have noted already in Paragraph~\ref{RTwist} that $\check{u}$ is invertible. By using
Lemma~\ref{QuasiHopf}, Lemma~\ref{QuasiTri}, and the counit-associator property, it follows directly from the definition that~\mbox{$\ea(u)=1$}, or alternatively~$\ea(\check{u})=1$ from its definition. The second property
of the Drinfel'd element is just one of the properties of~$\check{u}$ that follow directly from its construction in Paragraph~\ref{RTwist}. For the third property, recall that we have described the structure elements of the coopposite quasi-Hopf algebra in Paragraph~\ref{Coop}; in particular, we have seen there that the element~$F$, formed in~$A^{\cop}$, is just 
$(\sa^{-1} \o \sa^{-1})(F)$. On the other hand, we have explained in Paragraph~\ref{RTwist} how the coopposite coproduct arises by twisting the original coproduct with the help of the R-matrix, or alternatively with the help of its variant~$R'^{-1}$. As the two structures
were related via~$\hat{u}$ resp.~$\check{u}$, we get from
Proposition~\ref{ModAnti} that 
$$\hat{F} = (\hat{u} \o \hat{u}) (\sa^{-1} \o \sa^{-1})(F) \Delta^{\cop}(\hat{u}^{-1}) \qquad 
\check{F} = (\check{u} \o \check{u}) (\sa^{-1} \o \sa^{-1})(F) \Delta^{\cop}(\check{u}^{-1})$$
Because the Drinfel'd element is equal to~$\check{u}$, we focus on 
the second formula, and substitute for~$\check{F}$ the expression
from the beginning of this paragraph to get
$$(\sa \o \sa)(R) F R' = (u \o u) (\sa^{-1} \o \sa^{-1})(F) \Delta^{\cop}(u^{-1})$$
But we have $(\sa \o \sa)(R) F = F' R$ by Lemma~\ref{QuasiTri}, and therefore can use the second property of the
Drinfel'd element to rewrite the preceding equation as
$$F' R R' = (\sa \o \sa)(F) (u \o u) \Delta^{\cop}(u^{-1})$$
Interchanging tensor factors, this becomes
$F R' R = (\sa \o \sa)(F') (u \o u) \Delta(u^{-1})$,
which in turn implies 
$R'R \Delta(u) = F^{-1} (\sa \o \sa)(F') (u \o u)$.
But by quasi-cocom\-mutativity, we have $R'R \Delta(u) = \Delta(u) R'R$, and the third assertion follows.~\qed
\end{pf}

It must be emphasized that the preceding proposition is not new: The invertibility of~$u$, the first property and in particular the second property were proved by D.~Altsch\"uler and A.~Coste in \cite{AltCoste}, Sec.~3, p.~87f. The third property is stated there as well (cf.~Eq.~(4.21), p.~95), and the authors also propose a general strategy for its proof, of which they carry out the first step explicitly (cf.~Eq.~(4.20), p.~95), which however, as they say clearly, only works under the assumption that~$\alpha$ is invertible. 
The first complete, rigorous proof without this assumption was given
by D.~Bulacu and E.~Nauwelaerts in \cite{BulNauRib}, p.~668ff. As its Hopf-algebraic predecessor (cf.~\cite{M}, Thm.~10.1.13, p.~181f), it is based on a comparatively
involved computation, but has the advantage to deduce the result almost directly from the axioms.

\subsection[Ribbon quasi-Hopf algebras]{} \label{RibQuasiHopf}
We now use the machinery developed so far to study ribbon quasi-Hopf algebras. A quasitriangular quasi-Hopf algebra is called a ribbon quasi-Hopf algebra if it contains a ribbon element. This means the following:
\begin{defn}
A nonzero central element $v \in A$ is called a ribbon element if it satisfies
$$\da(v) = (R'R) (v \o v) \quad \text{and} \quad \sa(v) = v$$
\end{defn}

Let us clarify how this definition relates to the various competing
definitions of a ribbon quasi-Hopf algebra that we have already mentioned
in the introduction. We will prove below that it follows from our definition that a ribbon element is invertible. The definitions given in \cite{AltCoste}, \cite{BulNauRib}, and~\cite{BulPanOystTrace} all
work instead with the inverse element; our convention is the one used in \cite{Tur}, Sec.~XI.3.1, p.~500. As already pointed out in \cite{BulPanOystTrace}, Def.~2.3, p.~6106, it follows from the counitality property and Lemma~\ref{QuasiTri} that~$\ea(v)=1$; to see this, one just needs to apply $\ea \o \id$ to the first axiom in our definition above. This shows that, modulo the inversion, our definition matches with the definition in~\cite{BulPanOystTrace}, loc.~cit. 

A different definition was given by D.~Altsch\"uler and A.~Coste in \cite{AltCoste}, Par.~4.1, p.~89. As noted in \cite{BulNauRib}, Thm.~3.1,
p.~667, it follows from the formula for the coproduct of the Drinfel'd element, which we have just reconfirmed in Proposition~\ref{PropDrinf}.3, that the definition given by Altsch\"uler and Coste is equivalent to our definition and the additional requirement that $v^{-2} = u \sa(u)$. Furthermore, it was shown in \cite{BulPanOystTrace}, Prop.~5.5, p.~6119 that this property is automatically satisfied if~$\alpha$ is invertible. We will now show that this restriction is unnecessary. For preparation, we need the following lemma:
\begin{lemma} 
We have $v^2 \check{\alpha} = \hat{\alpha}$ 
and~$v^2 \hat{\beta} = \check{\beta}$.
\end{lemma}
\begin{pf}
Because the ribbon element is central and invariant under the antipode,
we have
\begin{align*}
v^2 \hat{\beta} & = \sum_{l} v^2 s_l \beta \sa(t_l) 
= \sum_{l} s_l v \beta \sa(t_l v) 
\end{align*}
The above definition also yields $R'^{-1} \da(v) = R (v \o v)$. Inserting this into the preceding formula, we get
\begin{align*}
v^2 \hat{\beta} &
= \sum_{l} \bar{t}_l v_\1 \beta \sa(\bar{s}_l v_\2)  
= \sum_{l} \bar{t}_l v_\1 \beta \sa(v_\2)  \sa(\bar{s}_l)  
= \sum_{l}  \bar{t}_l \beta \sa(\bar{s}_l) = \check{\beta}
\end{align*}
by the right antipode equation and the fact that $\ea(v)=1$, which we
already recorded above. This proves the second assertion. The proof of the first assertion is similar: Since 
$\da(v) R^{-1} = (v \o v) R' $, we have
\begin{align*}
v^2 \check{\alpha} = \sum_l \sa(v t_l) \alpha v s_l
= \sum_l \sa(v_\1 \bar{s}_l) \alpha v_\2 \bar{t}_l
= \sum_l \sa(\bar{s}_l) \alpha  \bar{t}_l = \hat{\alpha}
\end{align*}
by the left antipode equation.
\qed
\end{pf}

The proof of our main result is now almost immediate:
\begin{thm}
$v^{-2} = u \sa(u)$
\end{thm}
\begin{pf}
By construction, we have
$\sa^{-1}(\alpha) = \hat{u}^{-1} \hat{\alpha} = 
\check{u}^{-1} \check{\alpha}$.
Comparing this with the first assertion of the lemma, we see that 
$v^2 \check{\alpha} = \hat{\alpha} = \hat{u} \check{u}^{-1} \check{\alpha}$.
Now the duality axiom for the twisted quasi-Hopf algebra yields
$$\sum_i S(Z_i) \check{\alpha} Y_i \check{\beta} S(X_i) = 1$$
Because both $v^2$ and $\hat{u} \check{u}^{-1}$ are central, this implies
\begin{align*}
v^2 = \sum_i S(Z_i) v^2 \check{\alpha} Y_i \check{\beta} S(X_i) = 
\sum_i S(Z_i) \hat{u} \check{u}^{-1} \check{\alpha} Y_i \check{\beta} S(X_i)
= \hat{u} \check{u}^{-1}
\end{align*}
In view of Proposition~\ref{RTwist}, this means that
$v^2 = \sa^{-1}(u^{-1}) u^{-1}$. Inverting this, we get
$v^{-2} = u \sa^{-1}(u)$. But as~$u$ is
invariant under the square of the antipode by Proposition~\ref{PropDrinf}.2, this implies the assertion.
\qed
\end{pf}

\subsection[The opposite and coopposite quasi-Hopf algebra]{} \label{OpCoop}
In Paragraph~\ref{Coop}, we have described how to turn the coproduct into the coopposite coproduct. But we can also simultaneously turn the product into the opposite product. In this way, we arrive at the opposite and coopposite quasi-Hopf algebra~$A^{\op \cop}$, which is again a quasi-Hopf algebra with respect to the following structure elements: Its counit and antipode are unchanged, but its associator 
is $\sum_i Z_i \o Y_i \o X_i$, its evaluation element is~$\beta$, and its coevaluation element is~$\alpha$ (cf.~\cite{DrinfQuasiHopf}, Rem.~4, p.~1424; \cite{Kas}, Exerc.~XV.6.2, p.~381). Furthermore, if~$A$ was quasitriangular, then $A^{\op \cop}$ is still quasitriangular with respect to the same R-matrix. Therefore, its Drinfel'd element is
$$\tilde{u} := 
\sum_{i,l} \bar{Z}_i s_l \beta S(t_l) S(S(\bar{X}_i)\alpha\bar{Y}_i)$$
All the elements that we have introduced in Paragraph~\ref{RTwist} can also be formed in~$A^{\op \cop}$. But it turns out that we do not get any new elements in this way; rather these elements coincide with other elements formed in~$A$. For example, the element~$\hat{\alpha}$, if formed in~$A^{\op \cop}$, is equal to the original element~$\check{\beta}$ as formed in~$A$. The following table indicates which elements formed in~$A^{\op \cop}$ are equal to which elements formed in~$A$:
\begin{center}
\setlength{\extrarowheight}{2.7pt}
\begin{tabular}[t]{|l|c|c|c|c|c|c|}\hline 
In~$A^{\op \cop}$  & $\hat{\alpha}$ & $\hat{\beta}$ & $\check{\alpha}$ & $\check{\beta}$ & $\hat{u}$ & $\check{u}$ \\ \hline 
In~$A$ & $\check{\beta}$  & $\check{\alpha}$ & $\hat{\beta}$ & $\hat{\alpha}$ & $\check{u}^{-1}$ & $\hat{u}^{-1}$ \\ \hline 
\end{tabular} 
\end{center}
These correspondences can be applied to prove the following fact:
\begin{prop}
$u = \sa(\tilde{u})$
\end{prop}
\begin{pf}
By Proposition~\ref{RTwist}, we have $u=\check{u}$. In
$A^{\op \cop}$, this means~$\tilde{u} = \hat{u}^{-1}$.
But we have already seen in Proposition~\ref{RTwist} that
$\check{u} = \sa(\hat{u}^{-1})$.
\qed
\end{pf}

This result can also be proved by direct computation, which is quite tedious. However, there is another comparatively short proof: The result follows from Lemma~\ref{Twist}, because~$A^{\op \cop}$ is isomorphic to a twist of~$A$ by \cite{DrinfQuasiHopf}, Prop.~1.2, p.~1426.
Let us explain this in greater detail. The element~$\gamma$ introduced in Paragraph~\ref{QuasiHopf} satisfies 
$(\ea \o \id)(\gamma) = (\id \o \; \ea)(\gamma) = \ea(\alpha) \alpha$
by the counit-associator property. It then follows from the duality axiom that the element~$F$, which we have also defined there, satisfies
$(\ea \o \id)(F) = (\id \o \; \ea)(F) = \ea(\alpha) \A$, so that the element $T:= \ea(\beta) F$ satisfies the requirement
$(\ea \o \id)(T) = (\id \o \; \ea)(T) = \A$ imposed in Paragraph~\ref{Twist}; recall that $\ea(\alpha) \ea(\beta) = 1$ by Lemma~\ref{QuasiHopf}. As explained in \cite{DrinfQuasiHopf}, loc.~cit., the compatibility conditions stated in Proposition~\ref{QuasiHopf} now yield that the antipode, considered as a map from~$A^{\op\cop}$ to~$A_T$, is a quasi-bialgebra isomorphism.
However, it is not a quasi-Hopf algebra isomorphism; we rather have
$$\sa(\beta) = \ea(\beta)^2 \alpha_T \qquad \qquad 
\sa(\alpha) = \ea(\alpha)^2 \beta_T$$
as we see from \cite{BulNauRib}, Eq.~(2.14), p.~665 via a small correction. This means that the antipode becomes a quasi-Hopf algebra morphism if the evaluation element and the coevaluation element
of~$A_T$ are adjusted as indicated in Paragraph~\ref{ModAnti}, using the element $x := \ea(\beta)^2 \A$.

On the other hand, the compatibility between the antipode and the R-matrix stated in Lemma~\ref{QuasiTri} then yields that the antipode is in fact an isomorphism of quasitriangular quasi-Hopf algebras. It therefore maps the Drinfel'd element of~$A^{\op \cop}$ to the Drinfel'd element of~$A_T$, with the adjustments just indicated. By Proposition~\ref{QuasiTri}, this means in formulas that~$\sa(\tilde{u}) = x \sa(x^{-1}) u_T$.
However, we have $\sa(x) = x$ in our case, and therefore~$\sa(\tilde{u}) = u_T$. But $u_T=u$ by Lemma~\ref{Twist}, which completes the second derivation of our proposition above.

\end{document}